\documentclass{amsart}
\usepackage{amsrefs}
\title{Formulas generalizing Pappus and Desargues}
\author{Roger D. Maddux}
\address{Department of Mathematics\\ 396 Carver Hall\\ Iowa State
  University\\ Ames, Iowa 50011\\ U.S.A.}
\email{maddux@iastate.edu}
\date{17 November 2016}
\def\u{\mathbf u}
\def\v{\mathbf v}
\def\w{\mathbf w}
\def\x{\mathbf x}
\def\y{\mathbf y}
\def\z{\mathbf z}
\def\p{\mathbf p}
\def\q{\mathbf q}
\def\r{\mathbf r}
\def\0{\mathbf 0}

\let\bm=\bmatrix
\let\ebm=\endbmatrix
\let\vm=\vmatrix
\let\evm=\endvmatrix
\def\det#1#2#3{\mathsf{det}[{#1,\,}{#2,\,}{#3}]}
\def\ee#1#2#3{\bmatrix{#1}\\{#2}\\{#3}\endbmatrix}
\def\xyz#1#2#3{x_{#1}y_{#2}z_{#3}}

\def\wow#1#2#3#4{\det{#4}{#1}{#3}#2+\det{#3}{#2}{#4}#1}
\def\display#1{$$\text{#1}$$}
\def\<{\left<}
\def\>{\right>}
\def\({\left(}
\def\){\right)}
\def\1{\det\x\v\y}
\def\2{\det\x\w\z}
\def\3{\det\y\u\x}
\def\4{\det\y\w\z}
\def\5{\det\z\u\x}
\def\6{\det\z\v\y}
\def\reals{\mathbb F}

\def\solid{{\reals^3}}
\newtheorem{theorem}{Theorem}

\newtheorem{definition}{Definition}
\begin{document}
\allowdisplaybreaks
\begin{abstract}
  The Theorems of Pappus and Desargues are generalized by two special
  formulas that hold in the three-dimensional vector space over a
  field.
\end{abstract}
\maketitle
\section{Introduction}
The Theorems of Pappus and Desargues hold in the projective plane over
an arbitrary field $\reals$.  This projective plane can be constructed
in the three-dimensional vector space $\reals^3$ by taking the points
and lines to be the one-dimensional and two-dimensional subspaces of
$\reals^3$.  Two special formulas involving cross products and
determinants in $\reals^3$ have the Theorems of Pappus and Desargues
as corollaries.  These formulas were suspected and found in the spring
semester of 1997 as a result of teaching the second semester of
college geometry from \cite{MR854104}. The material for this paper was
extracted from the lecture notes for Math 436: Geometry \cite{M97}.
The lecture notes were based on the textbook \cite{MR854104}.
\section{The vector space $\mathbb F^3$}
Let $\reals$ be an arbitrary field. The elements of $\solid$ are
written as column vectors.
\begin{definition}
  For all $\x=\ee{x_1}{x_2}{x_3}\in\solid$,
  $\y=\ee{y_1}{y_2}{y_3}\in\solid$, and
  $\z=\ee{z_1}{z_2}{z_3}\in\solid$,
\begin{itemize}
\item the {\bf determinant} of $\x,\y,\z$ is
  \begin{align*}
    \det\x\y\z&=\vm x_1&y_1&z_1\\x_2&y_2&z_2\\x_3&y_3&z_3\evm
    \\&=\xyz123+\xyz312+\xyz231-\xyz132-\xyz321-\xyz213\in\reals,
  \end{align*}
\item the {\bf inner product} of $\x$ and $\y$ is
  \begin{align*}
    \<\x,\y\>=x_1y_1+x_2y_2+x_3y_3\in\reals,
  \end{align*}
\item the {\bf cross product} of $\x$ and $\y$ is
  \begin{align*}
    \x\times\y=\bm
    x_2y_3-x_3y_2\\x_3y_1-x_1y_3\\x_1y_2-x_2y_1\ebm\in\solid.
  \end{align*}
\end{itemize}
For arbitrary sets of vectors $X,Y\subseteq\reals^3$, $X+Y$ is their
complex sum: \display{$X+Y=\{\x+\y:\x\in X,\,\y\in Y\}$.}
\end{definition}
To prove the Pappus and Desargues Theorems we use the following pair
of specially designed formulas involving cross products and
determinants.
\begin{theorem}\label{th1}
  For any field $\reals$ and any $\u,\v,\w,\x,\y,\z\in\solid$,
  \begin{align}
    \label{P}\tag{P}
  &\det{{(\u\times\y)\times(\x\times\v)}}
          {{(\w\times\x)\times(\z\times\u)}}
          {{(\v\times\z)\times(\y\times\w)}}\\\notag
          &=\det\v\u\x\det\u\w\z\det\w\v\y\det\y\x\z \\\notag
          &\,+\det\x\y\v\det\z\x\u\det\y\z\w\det\u\w\v,\\
  \label{D}\tag{D}
  &\det{(\w\times\u)\times(\z\times\x)}
        {(\u\times\v)\times(\x\times\y)}
        {(\v\times\w)\times(\y\times\z)} \\\notag
        &=\det\x\y\z\det{\u\times\x}{\v\times\y}{\w\times\z}\det
        \u\v\w.
\end{align}
\end{theorem}
The proof of Theorem \ref{th1} could be left to the reader since it is
merely computational: for each formula, expand both sides in 18
variables (3 for each of 6 vectors) and see that the expansions are
the same.  Nevertheless, a proof of Theorem \ref{th1} using simpler
properties of cross products and determinants is included as an
appendix.
\section{The projective plane over $\mathbb F$.}
The projective plane over $\reals$ has points, lines, and an incidence
relation between them: a point may be on a line.
\begin{definition}
  The {\bf points} $X$, $Y$, \dots\ of the {\bf projective plane over}
  $\reals$ are the one-dimensional subspaces of the three-dimensional
  vector space $\reals^3$ over $\reals$, and its {\bf lines}
  $\lambda$, $\mu$, \dots\ are the two-dimensional subspaces of
  $\reals^3$.  A point is {\bf on} a line if it is contained in that
  line: \display{$X$ is {\bf on} $\lambda$ iff $X\subseteq\lambda$.}
\end{definition}
Every non-zero vector in $\solid$ determines both a point and a line.
Suppose $0\neq\x,\y\in\solid$. Let \display{$X=\{r\x:r\in\reals\}$ and
  $\lambda=\{\z:\<\z,\y\>=0\}$.}  Then $X$ is the 1-dimensional
subspace generated by $\x$ and $\lambda$ is the 2-dimensional subspace
of vectors perpendicular to $\y$, so $X$ is a point and $\lambda$ is a
line in the projective plane over $\reals$.  The vectors $\x$ and $\y$
are called (sets of) {\bf homogeneous coordinates} for $X$ and
$\lambda$, respectively.  The incidence relation between a point and a
line holds if the inner product of their homogeneous coordinates is
$0$ (their homogeneous coordinates are perpendicular): \display{$X$ is
  on $\lambda$ iff $\<\x,\y\>=0$.}  Suppose $X$ and $Y$ are distinct
points.  Then there are non-zero vectors $\x,\y\in\solid$ such
that
\begin{equation*}
  X=\{r\x:r\in\reals\}\text{ and }Y=\{s\y:s\in\reals\}.
\end{equation*}
The unique line containing both $X$ and $Y$ is the 2-dimensional
subspace generated by $\x$ and $\y$:
\begin{equation*}
  X+Y=\{r\x+s\y:r,s\in\reals\}.
\end{equation*}
The cross product $\x\times\y$ is perpendicular to both $\x$ and $\y$,
so any vector perpendicular to $\x\times\y$ will be in the unique line
$X+Y$ that contains $X$ and $Y$:
\display{$X+Y=\{\z:\<\z,\x\times\y\>=0\}$.}  Thus points with
homogeneous coordinates $\x$ and $\y$ determine a line with
homogeneous coordinates $\x\times\y$.  Suppose two lines, $\lambda$
and $\mu$, have homogeneous coordinates $\x,\y\in\reals$,
respectively:
\begin{equation*}
  \lambda=\{\z:\<\z,\x\>=0\},\qquad\mu=\{\z:\<\z,\y\>=0\}.
\end{equation*}
Note that $\x$ and $\y$ are linearly independent, since otherwise
$\lambda=\mu$ (and there is only one line, not two, as postulated).
Then $\lambda\cap\mu$ is a 1-dimensional subspace of $\reals^3$ (a
point) generated by any non-zero vector perpendicular to both $\x$ and
$\y$, such as $\x\times\y$:
\begin{equation*}
  \lambda\cap\mu=\{t(\x\times\y):t\in\reals\}.
\end{equation*}
Thus, two lines with homogeneous coordinates $\x$ and $\y$ meet at a
point with homogeneous coordinates $\x\times\y$.

Suppose $\0\neq\x,\y,\z\in\reals^3$.  Let $X$, $Y$, and $Z$ be the
points, and $\lambda$, $\mu$, and $\nu$ the lines, with homogeneous
coordinates $\x$, $\y$ and $\z$, respectively:
\begin{align*}
X&=\{r\x:r\in\reals\},&\quad 
Y&=\{r\y:r\in\reals\},&\quad 
Z&=\{r\z:r\in\reals\},\\
\lambda&=\{\u:\<\u,\x\>=0\},&\quad
\mu&=\{\u:\<\u,\y\>=0\},&\quad
\nu&=\{\u:\<\u,\z\>=0\}.
\end{align*}
By definition, $X,Y,Z$ are {\bf collinear} iff $X+Y=X+Z=Y+Z$, and
$\lambda,\mu,\nu$ are {\bf concurrent} iff $\lambda\cap\mu\cap\nu$ is
a point. Concurrency and collinearity are equivalent to the
homogeneous coordinates having determinant zero:
\begin{equation*}
  X,Y,Z\text{ are collinear} \quad\text{iff}\quad \det\x\y\z=0
  \quad\text{iff}\quad \lambda,\mu,\nu\text{ are concurrent.}
\end{equation*}
\section{The theorems of Pappus and Desargues}
It is convenient to state the Theorems of Pappus and Desargues
together because between them they deal with all fifteen lines passing
through six points (fifteen is six taken two at a time). Twelve of
these lines are intersected in pairs to produce six more points, three
for Pappus and the other three for Desargues, while the remaining
three lines are involved in the conclusion of Desargues.
\begin{theorem}[Pappus-Desargues]\label{th2}
  Let $U,V,W,X,Y,Z$ be six points in the projective plane over
  $\reals$. Define six more points
\begin{align*}
  O&=(V+Z)\cap(Y+W),&R&=(V+W)\cap(Y+Z),\\
  P&=(W+X)\cap(Z+U),&S&=(W+U)\cap(Z+X),\\
  Q&=(U+Y)\cap(X+V),&T&=(U+V)\cap(X+Y).
\end{align*}
\begin{enumerate}
\item[] \hskip-.4in{\rm(\bf Pappus)} If $U,V,W$ are collinear and
  $X,Y,Z$ are collinear then $O,P,Q$ are collinear.
\item[] \hskip-.4in{\rm(\bf Desargues)} If $U,V,W$ are not collinear
  and $X,Y,Z$ are not collinear then the lines $U+X$, $V+Y$, $W+Z$ are
  concurrent iff the points $R,S,T$ are collinear.
\end{enumerate}
\end{theorem}
\proof Choose non-zero vectors $\u,\v,\w,\x,\y,\z\in\solid$ that are
homogeneous coordinates for the points $U,V,W,X,Y,Z$, respectively:
$U=\{r\u:r\in\reals\}$, {\it etc}.  From the definitions of the
additional six points we obtain homogeneous coordinates for them as
cross products of cross products:
\begin{align*}
    O&=\{r\((\v\times\z)\times(\y\times\w)\):r\in\reals\},\\
    P&=\{r\((\w\times\x)\times(\z\times\u)\):r\in\reals\},\\
    Q&=\{r\((\u\times\y)\times(\x\times\v)\):r\in\reals\},\\
    R&=\{r\((\v\times\w)\times(\y\times\z)\):r\in\reals\},\\
    S&=\{r\((\w\times\u)\times(\z\times\x)\):r\in\reals\},\\
    T&=\{r\((\u\times\v)\times(\x\times\y)\):r\in\reals\}.
\end{align*}
For Pappus's Theorem, assume $U,V,W$ are collinear and $X,Y,Z$ are
collinear.  Then
\begin{equation*}
  \det\y\x\z=0=\det\u\w\v,
\end{equation*}
so by formula \eqref{P} we get
\begin{align*}
  &\det{(\u\times\y)\times(\x\times\v)}{(\w\times\x)\times(\z\times\u)}
  {(\v\times\z)\times(\y\times\w)}\\
  &=\nonumber\det\v\u\x\det\u\w\z\det\w\v\y\underbrace{\det\y\x\z}_0\\
  &\,+\nonumber\det\x\y\v\det\z\x\u\det\y\z\w\underbrace{\det\u\w\v}_0=0,
\end{align*}
hence $O,P,Q$ are collinear. For Desargues's Theorem, assume that
$U,V,W$ are not collinear and $X,Y,Z$ are not collinear. Hence
\begin{equation*}
  \det\u\v\w\neq0\neq\det\x\y\z.
\end{equation*}
The conclusion of Desargues's Theorem is that the following two
statements are equivalent (the first one says $U+X$, $V+Y$, $W+Z$ are
concurrent, the second says $R,S,T$ are collinear):
\begin{align}
  0&=\det{\u\times\x}{\v\times\y}{\w\times\z},\label{1st}\\
  0&=\det{(\w\times\u)\times(\z\times\x)}{(\u\times\v)\times(\x\times\y)}
  {(\v\times\w)\times(\y\times\z)}.\label{2nd}
\end{align}
This equivalence follows immediately from formula \eqref{D}, written
here with the non-collinearity assumptions included:
\begin{align*}
  &\det{(\w\times\u)\times(\z\times\x)}{(\u\times\v)\times(\x\times\y)}
  {(\v\times\w)\times(\y\times\z)}\\
  &=\underbrace{\det\x\y\z}_{\neq0}\det{\u\times\x}{\v\times\y}{\w\times\z}
  \underbrace{\det\u\v\w}_{\neq0}.
\end{align*}
\endproof Notice that \eqref{1st} implies \eqref{2nd} without the
non-collinearity assumptions but to get \eqref{1st} from \eqref{2nd}
requires knowing that if a product of three numbers is zero and two of
them are not zero then the third one must be zero.  Here the
non-collinearity assumptions are needed.  For example, if $X,Y,Z$ are
collinear then $R,S,T$ are also collinear because they lie on the same
line as $X,Y,Z$, but the lines $U+X$, $V+Y$, and $W+Z$ are free to not
concur.

With or without any of the assumptions of Pappus or Desargues, the
formulas \eqref{P} and \eqref{D} express explicit numerical
relationships holding among the six points and fifteen lines
connecting them.  The Theorems of Desargues and Pappus just deal with
cases in which some of the determinants in the formulas are zero. But
the formulas \eqref{P} and \eqref{D} hold all the time. They are in
this sense strict generalizations of the Theorems of Pappus and
Desargues for the projective planes arising from fields.  For example,
given six arbitrary points $U,V,W,X,Y,Z$ in the projective plane over
$\reals$, we can create four triples of points, one of which is
$O,P,Q$ from Theorem \ref{th2}, and the other three are
\begin{align*}
  &(Z+V)\cap(U+W),&&(Y+U)\cap(X+Z),&&(W+X)\cap(V+Y),\\
  &(Y+U)\cap(W+V),&&(X+W)\cap(Z+Y),&&(V+Z)\cap(U+X),\\
  &(Z+W)\cap(V+X),&&(U+V)\cap(Y+Z),&&(X+Y)\cap(W+U),
\end{align*}
with the property that homogeneous coordinates for these triples all
have the same determinant, by \eqref{P}.  Therefore either all four
triples are collinear or none of them are collinear. Pappus's Theorem
says only that if $U,V,W$ and $X,Y,Z$ are collinear then $O,P,Q$ and
the other three triples are collinear.
\section*{Appendix: Proof of Theorem \ref{th1}}
Here are the properties of determinants and cross products used in the
proofs of formulas \eqref{P} and \eqref{D}.
\begin{align}
  \intertext{Cross products are perpendicular to their factors:}
  \label{3}&\<\x\times\y,\x\>=\<\x\times\y,\y\>=0 \intertext{The scalar
    triple product (box product) formula:}
  \label{5}&\<\x\times\y,\z\> =\<\x,\y\times\z\>=\det\x\y\z
  \intertext{Determinants with proportional inputs are zero:}
  \label{6}&\det\x\y{r\x}=0 \intertext{Determinants are invariant under
    cyclic permutations:}
  \label{7}&\det\x\y\z=\det\y\z\x=\det\z\x\y \intertext{Switching
    determinant inputs changes sign:}
  \label{8}&\det\x\y\z=-\det\y\x\z=-\det\x\z\y=-\det\z\y\x
  \intertext{Determinants are unchanged by transposition:}
  \label{9}&\det\x\y\z =\vm x_1&y_1&z_1\\x_2&y_2&z_2\\x_3&y_3&z_3\evm
  =\vm x_1&x_2&x_3\\y_1&y_2&y_3\\z_1&z_2&z_3\evm=\mathsf{det}\ee\x\y\z
  \intertext{Scalars move in and out of determinants:}
  \label{10}&r\det\x\y\z =\det{r\x}\y\z=\det\x{r\y}\z=\det\x\y{r\z}
  \intertext{Determinants distribute over vector addition:}
  \label{11}&\det{\w+\x}\y\z =\det\w\y\z+\det\x\y\z \intertext{A
    product of determinants is the determinant of inner products:}
  \label{12}&\mathsf{det}\ee\u\v\w\det\x\y\z=\vm\<\u,\x\>&\<\u,\y\>&\<\u,\z\>\\
  \<\v,\x\>&\<\v,\y\>&\<\v,\z\>\\\<\w,\x\>&\<\w,\y\>&\<\w,\z\>\evm
  \intertext{The vector triple product formula:}
  \label{13}&(\x\times\y)\times\z =\<\x,\z\>\y-\<\y,\z\>\x
  \intertext{The vector quadruple product formula:}
  \label{15}&(\w\times\x)\times(\y\times\z)=\det\z\w\y\x+\det\y\x\z\w
\end{align}
Proof of the vector quadruple product formula:
\begin{align*}
  (\w\times\x)\times(\y\times\z)
  &=\<\w,\y\times\z\>\x-\<\x,\y\times\z\>\w&&\text{by \eqref{13}}\\
  &=\det\w\y\z\x-\det\x\y\z\w&&\text{by \eqref{5}}\\
  &=\det\z\w\y\x+\det\y\x\z\w&&\text{by \eqref{7}, \eqref{8}}
\end{align*}
Proof of formula \eqref{P}. Define three vectors and use the vector
quadruple product formula \eqref{15} to rewrite them.
\begin{align*}
  \p&={(\u\times\y)\times(\x\times\v)}=\wow\u\y\x\v,\\
  \q&={(\w\times\x)\times(\z\times\u)}=\wow\w\x\z\u,\\
  \r&={(\v\times\z)\times(\y\times\w)}=\wow\v\z\y\w.
\end{align*} 
Then $\det\p\q\r$, the left side of \eqref{P}, is the sum of eight
determinants by the distributive law \eqref{11}.  Move all the scalars
(the determinants) out by \eqref{10}, obtaining another sum of eight
terms. Six of these cancel with others as indicated by the notations
to the right of the affected terms. What remains are the two terms on
the right side of \eqref{P}.
\begin{align*}
  \det\p\q\r &=\text{\sf det}
  \Big[\det\v\u\x\y+\det\x\y\v\u,\\
  &\kern34pt\det\u\w\z\x+\det\z\x\u\w,\\
  &\kern34pt\det\w\v\y\z+\det\y\z\w\v \Big]\\
  &=\det{\det\v\u\x\y}{\det\u\w\z\x}{\det\w\v\y\z}\\
  &\quad+\det{\det\v\u\x\y}{\det\u\w\z\x}{\det\y\z\w\v}\\
  &\quad+\det{\det\v\u\x\y}{\det\z\x\u\w}{\det\w\v\y\z}\\
  &\quad+\det{\det\v\u\x\y}{\det\z\x\u\w}{\det\y\z\w\v}\\
  &\quad+\det{\det\x\y\v\u}{\det\u\w\z\x}{\det\w\v\y\z}\\
  &\quad+\det{\det\x\y\v\u}{\det\u\w\z\x}{\det\y\z\w\v}\\
  &\quad+\det{\det\x\y\v\u}{\det\z\x\u\w}{\det\w\v\y\z}\\
  &\quad+\det{\det\x\y\v\u}{\det\z\x\u\w}{\det\y\z\w\v}\\
  &=\det\v\u\x\det\u\w\z\det\w\v\y\det\y\x\z\\
  &\quad+ \det\v\u\x\det\u\w\z\det\y\z\w\det\y\x\v\quad=r_1=-r_5
  \text{ by \eqref{7}, \eqref{8}}\\
  &\quad+ \det\v\u\x\det\z\x\u\det\w\v\y\det\y\w\z\quad=r_2=-r_3
  \text{ by \eqref{7}, \eqref{8}}\\
  &\quad+ \det\v\u\x\det\z\x\u\det\y\z\w\det\y\w\v\quad=r_3\\
  &\quad+ \det\x\y\v\det\u\w\z\det\w\v\y\det\u\x\z\quad=r_4=-r_6
  \text{ by \eqref{8}}\\
  &\quad+ \det\x\y\v\det\u\w\z\det\y\z\w\det\u\x\v\quad=r_5 \\
  &\quad+ \det\x\y\v\det\z\x\u\det\w\v\y\det\u\w\z\quad=r_6 \\
  &\quad+ \det\x\y\v\det\z\x\u\det\y\z\w\det\u\w\v \\
  &=\det\v\u\x\det\u\w\z\det\w\v\y\det\y\x\z \\
  &\quad+ \det\x\y\v\det\z\x\u\det\y\z\w\det\u\w\v.
\end{align*}
Proof of formula \eqref{D}. Define three vectors and use the vector
quadruple product formula \eqref{15}.
\begin{align*}
  \p=(\w\times\u)\times(\z\times\x)&={\2\u+\5\w},\\
  \q=(\u\times\v)\times(\x\times\y)&={\3\v+\1\u},\\
  \r=(\v\times\w)\times(\y\times\z)&={\6\w+\4\v}.
\end{align*}
Then $\det\p\q\r$, the left side of  \eqref{D}, is the sum of eight
determinants by \eqref{11}, six of which are zero by \eqref{6} because
two of the argument vectors are proportional.
\begin{align*}
  \det\p\q\r&=\begin{aligned}\text{\sf det}[&{\2\u+\5\w},\\
    &{\3\v+\1\u},\\
    &{\6\w+\4\v}]\end{aligned}\\
  &=\det{\2\u}{\3\v}{\6\w}\\
  &\quad+\det{\2\u}{\3\v}{\4\v}&&\text{$=0$ by \eqref{6}}\\
  &\quad+\det{\2\u}{\1\u}{\6\w}&&\text{$=0$ by \eqref{6}}\\
  &\quad+\det{\2\u}{\1\u}{\4\v}&&\text{$=0$ by \eqref{6}}\\
  &\quad+\det{\5\w}{\3\v}{\6\w}&&\text{$=0$ by \eqref{6}}\\
  &\quad+\det{\5\w}{\3\v}{\4\v}&&\text{$=0$ by \eqref{6}}\\
  &\quad+\det{\5\w}{\1\u}{\6\w}&&\text{$=0$ by \eqref{6}}\\
  &\quad+\det{\5\w}{\1\u}{\4\v}\\
  &=\2\3\6\det\u\v\w\\
  &\quad+\5\1\4\det\w\u\v&&\text{by \eqref{10}}\\
  &=\Big(\2\3\6\\
  &\quad+\5\1\4\Big)\det\u\v\w&&\text{by \eqref{7}}\\
  &=\vm0&\1&\2\\\3&0&\4\\\5&\6&0\evm\det\u\v\w
  &&\text{by definition of $\mathsf{det}$}\\
  &=\vm\<\x,\u\times\x\>&\<\x,\v\times\y\>&\<\x,\w\times\z\>\\
  \<\y,\u\times\x\>&\<\y,\v\times\y\>&\<\y,\w\times\z\>\\
  \<\z,\u\times\x\>&\<\z,\v\times\y\>&\<\z,\w\times\z\>\evm\det\u\v\w
  &&\text{by \eqref{5} and \eqref{3}}\\
  &=\text{\sf det}\ee\x\y\z\det{\u\times\x}{\v\times\y}{\w\times\z}\det\u\v\w
  &&\text{by \eqref{12}}\\
  &=\det\x\y\z\det{\u\times\x}{\v\times\y}{\w\times\z}\det\u\v\w
  &&\text{by \eqref{9}}
\end{align*}

\end{document}